\documentclass[12pt,runningheads]{article}
{\begin{Sbox}\begin{minipage}}%
{\end{minipage}\end{Sbox}\fbox{\TheSbox}}%
\usepackage{marvosym}
\usepackage[psamsfonts]{amssymb}
\usepackage{amsmath, amsthm, anysize, enumerate, color, url}
\usepackage{graphicx}
\usepackage{epstopdf}
\usepackage{epsfig}
\usepackage{subfigure}
\usepackage{setspace}

\usepackage{algorithm}
\usepackage{algpseudocode}

\usepackage{cite}
\usepackage[authoryear]{natbib}

\usepackage{threeparttable}

\newtheorem{theorem}{Theorem} 
\newtheorem{corollary}{Corollary}

\usepackage[colorlinks, breaklinks,
                   linkcolor=red,
                  anchorcolor=blue,
                  citecolor=blue
                  ]{hyperref}

\usepackage{breakurl}

\theoremstyle{definition}
\newtheorem{remark}{Remark}  
\newtheorem{example}{Example}
\newtheorem{condition}{Condition}

\newcommand{\E}{\mathbb{E}}

\newcommand{\R}{\mathbb{R}}

\renewcommand{\P}{\mathbb{P}}

\usepackage{setspace}
\usepackage{bm}

\usepackage{amsmath}
\usepackage{amssymb}
\usepackage{titletoc}
\usepackage{amsfonts,latexsym,amsthm,amsxtra,bbm}
\allowdisplaybreaks
\usepackage{graphicx}
\usepackage{epstopdf}
\usepackage{chngcntr}
\usepackage{booktabs}
\usepackage[font=small]{caption}
\usepackage{url}
\usepackage{tabularx} 
\usepackage{siunitx} 
\counterwithin{figure}{section}

\allowdisplaybreaks

\begin{document}

\title{Inference in covariate-adjusted bipartite network models} 

\author{
Zuhui Wu\thanks{Department of Statistics, Central China Normal University, Wuhan, 430079, China.}
\hspace{4mm}
Qiuping Wang\thanks{School of Mathematics and Statistics, Zhaoqing University,Zhaoqing,526000,China.
\texttt{Emails:} qp.wang@mails.ccnu.edu.cn.}
\hspace{4mm}
Ting Yan\thanks{Department of Statistics, Central China Normal University, Wuhan, 430079, China.}
}

\date{}

\maketitle

\begin{abstract}
In this paper, we introduce a general model for
jointly modelling the nodal heterogeneity and  covariates
in weighted or unweighted bipartite networks, which contains two different types of nodes.
The model has a degree heterogeneity parameter for each node and
 a fixed-dimensional regression coefficient for the covariates.
We use the method of moments to estimate the unknown parameters.
When the model belongs to the exponential family of distributions,
the moment estimator is identical to the maximum likelihood estimator.
We show the uniform consistency of the moment estimator, when the number of actors and
the number of events both go to infinity under some conditions. Further, we derive an asymptotic
representation of the moment estimator, which leads to their asymptotic
normal distributions under some conditions.
We present two applications to illustrate the unified results.
Numerical simulations and a real-data analysis demonstrates our theoretical findings.


\end{abstract}

\vskip 5pt

\section{Introduction}

A bipartite network is a kind of graph,
whose nodes are divided into two distinct sets, with edges only connecting nodes from different sets. Typical examples include
assigning tasks to workers or matching students to schools in matching problems,
recommending items to users in recommendation systems,
linking users to events in social networks,
business relationship between buyers and sellers in economic markets,
 connections between genes and the proteins in biological networks,
 affiliations between terms and documents in information retrieval.
Sometimes, a bipartite network is called an affiliation network, a bipartite graph or bigraph.
One of tasks in the analysis of bipartite networks is to understand
the generative mechanism of link formations, which can help to address
a diverse set of societal and economic issues including
information dissemination, item or product recommendation, efficiency of communication,
among many others \citep[e.g.][]{skvoretz1999logit,kunegis2010link,todd2015applying,kaya2020hotel,  liew2023methodology}

The nodal heterogeneity is commonly observed in real-world networks, which
describes variation in the abilities of different nodes
to participate in network connection. It is measured by the degrees of nodes,
which count the number of links connecting other nodes.
In bipartite networks, there are two sets of degrees for two sets of different nodes.
We use  ``event" and ``actor" to denote these two kinds of nodes hereafter.
Some ``events" in one node set may have links to many ``actors" in the other node set,
 while a large number of other events have relatively fewer edges.
\cite{Wang2022JMLR} proposed a bipartite $\beta$-model to model the nodal heterogeneity
in bipartite networks, which is an exponential family distribution with the node degrees
as the sufficient statistics, with each node being assigned one node-specific parameter.
It extends the $\beta$-model \citep{Chatterjee:Diaconis:Sly:2011} for undirected graphs
to bipartite graphs, which can be date back to an earlier $p_1$ model for directed networks \citep{Holland:Leinhardt:1981}.
\cite{fan2022asymptotic} further
established the consistency and asymptotic normality
of the moment estimator when the numbers of ``actors" and ``events" both go to infinity
in a general version of the bipartite $\beta$-model.

The attributes of nodes also have a significant influence on the network generation process.
The actors  are more likely to connect with the events whose attributes are
in coordinate with those of actors,
than other events with inconsistent attributes. For instance, in a social event recommendation network, if a user (actor) has the attributes of ``loves outdoor activities", ``enjoys hiking", and ``is interested in environmental protection", this user is much more likely to establish connections with outdoor hiking events, forest cleanup volunteer activities, or eco-friendly camping events, rather than indoor e-sports competitions, painting exhibitions, or indoor concert events.

Several models have been proposed to quantify the influence of
the nodal heterogeneity and attributes of nodes on the network formation,
including a covariate-adjusted $\beta$-model for undirected networks  \citep{graham2017econometric}, a covariate-adjusted $p_0$ model for directed networks
\citep{yan2019} and a probit-type regression model \citep{dzemski2019empirical},
where the covariates are associated with edges and could form according the attributes of nodes in a paired manner.
In these models, an edge forms when the surplus made up of the effects of both nodes and covariates,
is over the latent random noise.
Specifically, \cite{graham2017econometric} and \cite{yan2019} assumed a logistic distribution
on the latent random noise,
while \cite{dzemski2019empirical} assumed a normal distribution.
Further, \cite{wang2023asymptotic}
establish a unified framework for
 a general model that jointly characterizes degree heterogeneity and
homophily in weighted, undirected networks and establish consistency and
asymptotic normality of the moment estimator.
However, a unified analysis for bipartite networks with covariates have not been explored.

In this paper, we introduce a general model for
jointly modelling the nodal heterogeneity and  covariates
in weighted or unweighted bipartite networks.
Assuming there are $m$ actors and $n$ events, the model contains
a set of degree parameters $\{\alpha_i\}_{i=1}^m$ for actors,
a set of degree parameters $\{\beta_j\}_{i=1}^n$ for events
and a fixed-dimensional regression coefficient $\gamma$ for the covariates.
The degree parameters measure the node heterogeneity while
the regression coefficient measure the effects of covariates.
We use the method of moments to estimate the unknown parameters, rather than
the likelihood approaches \citep{graham2017econometric,yan2019}.
As discussed later in the paper,
the likelihood equations under some common distributions such as the probit distribution
are not convenience to analyze the estimator in a unified manner.
Another reason is that the moment estimation can be extended into dependent edges,
although this paper exclusively focuses on independent edge generation.
When the model belongs to the exponential family of distributions,
the moment estimator is identical to the maximum likelihood estimator.
By developing a two-stage Newton method that first finds an error bound for $\hat{\alpha}_\gamma$
and $\hat{\beta}_\gamma$ with a fixed  $\gamma$ via establishing the convergence rate of the Newton iterative
sequence and then derives the error of $\hat{\gamma}_{\alpha,\beta}$ with fixed $\alpha$ and $\beta$,
we show the uniform consistency of the moment estimator, when the number of actors and
the number of events both go to infinity. Further, we derive an asymptotic
representation of the moment estimator, which leads to their asymptotic
normal distributions under some conditions.
We present two applications to illustrate the unified results.
Numerical simulations and a real-data analysis demonstrates our theoretical findings.

The rest of this paper is organized as follows. We first set up some notation; Section 2 describes
our proposed model, devises estimation methods and inference procedures, and establishes theoretical justifications; in Section 3, we conduct comprehensive simulation studies to address
the question of selecting the tuning parameter, assess the performance of our method from different aspects, and validate our theoretical prediction; in Section 4, we apply our method to
a real-world dataset and interpret the results; we conclude our paper with some discussion in
Section 5.

\section{A covariate-adjusted bipartite network model}
\label{sec-model}

Recall that we consider
a bipartite network $\mathcal{G}(m,n)$ consisting of $m$ actors labelled as $\{1, \ldots,m\}$ and $n$ events labelled as $\{1,\ldots,n\}$.
In real-world bipartite networks,
the number of actors is usually larger than the number of events.
Without loss of generality, we assume $m \geq n$ throughout the paper.
(The analysis is similar when $m<n$.)
We use $x_{ij}$ to denote the edge weight between actor $i$ and event $j$
and $X=(x_{ij})_{m\times n}$ to denote the adjacency matrix of $\mathcal{G}_{m,n}$.
The edge weight may be binary values $\{0,1\}$, nonnegative integer values
$\{0,1,\ldots\}$ or continuous nonnegative values $\mathbb{R}^+$.

Let $d=(d_1,\ldots,d_m)^\top$ and $b=(b_1,\ldots,b_n)^\top$ be the respective degree sequences of actors and events,
where $d_i=\sum_{j=1}^{n}x_{ij}$ and $b_j=\sum_{i=1}^{m}x_{ij}$.
Further, we collect the covariates associated with each edge.
We use $z_{ij}\in\mathbb{R}^p$ to denote the covariate of the edge between actor $i$ and event $j$. We assume that the dimension $p$ is fixed.

We use the degree parameters $\{\alpha_i\}_{i=1}^m$ to measure
the nodal heterogeneity for $m$ actors to participate in network connections.
For $n$ events, the corresponding node heterogeneity parameters are $\{\beta_j\}_{j=1}^n$.
The degree parameter $\alpha_i$ of actor $i$ measures the actor's activity level,
while  the degree parameter $\beta_j$ of event $j$  measures the event's popularity.
Recall that we use a $p$-dimensional parameter vector $\gamma$
to measure the effects of covariates.
Our goal is to jointly model nodal heterogeneity and the effects of covariates.
We formulate the covariate-adjusted bipartite network model as follows:
\begin{equation}
\label{model}
    x_{ij}|\{z_{ij}, \alpha, \beta, \gamma\}
    \sim
    f\big( x_{ij} \big| \alpha_i + \beta_j + z_{ij}^\top \gamma \big),
\end{equation}
where $f(\cdot)$ is a known probability mass or density function.
Here, the degree parameters and the covariate term enter into the model in an additive manner.
The parameter $\gamma\in \R^p$ is known as a regression coefficient of covariates.
Conditional on all covariates $\{z_{ij}\}_{i,j}$, we assume that all edges are independent.
Two running examples for illustrating the model are given below.

\begin{example}(Binary weight)
\label{example-a}
Consider $x_{ij}\in \{0,1\}$, indicating whether there is an edge between actor $i$ and
event $j$.
Then, the probability distribution of $x_{ij}$ is
\[
\P( x_{ij}= a ) = [F( \alpha_i+\beta_j + z_{ij}^\top \gamma )]^{a} (1-F(\alpha_i+\beta_j+z_{ij}^\top \gamma ) )^{1-a},~~a=0, 1.
\]
where $F$ is the cumulative distribution function of $f(\cdot)$.
Two common examples for $F(\cdot)$ are the logistic distribution: $F(x)=e^x(1+e^x)^{-1}$ and probit distribution: $F(x)=\Phi(x)$.
Here, $\Phi(x)$ is the cumulative distribution function for the standard normal random variable.
\end{example}

\begin{example}(Infinite discrete weight)
\label{example-b}
Consider $a_{ij}\in \{0, 1, \ldots\}$.
We assume $x_{ij}$ is from a Poisson distribution with mean $\lambda= \exp( \beta_i+\beta_j + z_{ij}^\top \gamma)$, i.e.,
\[
\log \P(a_{ij}=a) = a(\alpha_i+\beta_j + z_{ij}^\top \gamma) - \exp ( \alpha_i+\beta_j + z_{ij}^\top \gamma) - \log a!.
\]
\end{example}

It is easy to see that adding a constant to $\alpha_i$ and subtracting a same constant from $\beta_j$ leads to the same probability density function.
Therefore, the model is not identifiable without any restriction. Following Yan et al. (2016), we
set
\begin{equation}
\beta_n = 0,
\end{equation}
for the model identifiability.

\section{Moment estimation}
\label{sec-moment}

We use the method of moments to estimate the unknown model parameters.
Let $\mu(\cdot)$ denote the expectation of $f(\cdot)$.
For notational convenience, we define
\begin{equation}
\label{eq-definition-nota}
\begin{array}{c}
\theta:=(\alpha_{1},\ldots,\alpha_{m},\beta_{1},\ldots,\beta_{n-1}),
\\
\pi_{ij}:= \alpha_i + \beta_j  + z_{ij}^\top \gamma,
\\
\mu_{ij}(\theta,\gamma):=\mu(\alpha_{i}+\beta_{j} + z_{ij}^\top \gamma).
\end{array}
\end{equation}
Since the distribution of $x_{ij}$ depends only on $\pi_{ij}$, we could write
\[
\mathbb{E}(x_{ij})=\mu(\pi_{ij})=\mu( \alpha_i + \beta_j + z_{ij}^\top \gamma ).
\]
The moment equations can be written as
\begin{equation}
\label{eq-moment}
\begin{array}{rcl}
d_{i}  &  = &  \sum_{k=1}^{n}\mu(\alpha_{i}+\beta_{k}+z_{ik}^\top \gamma),  \quad  i=1,  \cdots, m,
\\
b_{j}  &  = &  \sum_{k=1}^{m}\mu(\alpha_{k}+\beta_{j}+z_{kj}^\top \gamma),  \quad  j=1,  \cdots, n-1,
\\
\sum_{i=1}^m\sum_{j=1}^n z_{ij}x_{ij} &  =
& \sum_{i=1}^m\sum_{j=1}^n z_{ij}\mu( \alpha_i + \beta_j + z_{ij}^\top \gamma ).
\end{array}
\end{equation}

The solution to the moment equations in \eqref{eq-definition-nota},
 i.e., the moment estimator, is denoted as $(\widehat\theta,\widehat\gamma)$, where $\widehat\theta=(\widehat{\alpha}_{1},\cdots,\widehat{\alpha}_{m},\widehat{\beta}_{1},\cdots,\widehat{\beta}_{n-1})$
and $\widehat{\beta}_{n}=0$.
We use the superscript ``*" to denote the true parameter values, i.e., $(\theta^*,\gamma^*)$
is the true parameter of $(\theta,\gamma)$.

We explain why we do not use the maximum likelihood equation here.
Take the probit distribution  in Example \ref{example-a} for example.
The maximum likelihood equation for the parameter $\alpha_i$ is
\[
\sum_{j\neq i} \left[ \frac{ a_{ij} \phi(\alpha_i + \beta_j + z_{ij}^\top \gamma) }{ \Phi(\alpha_i + \beta_j + z_{ij}^\top \gamma) }
- \frac{(1-a_{ij})\phi(\alpha_i + \beta_j + z_{ij}^\top \gamma) }{ 1-\Phi(\alpha_i + \beta_j + z_{ij}^\top \gamma) } \right] =0,
\]
where $\phi(\cdot)$ is the density function of the standard normality.
As we can see, the maximum likelihood equation
is far from the moment equation, which has  a unified form in \eqref{eq-moment}.
Also, the analysis of this equation is more challenging since the random variables
are involved with the Jacobian matrix of the function made up of the formula of the left side of the equation.

\section{Asymptotic properties}
\label{sec-asymptotic}

In this section, we present the consistency and asymptotic normality of the moment estimator.

\subsection{Preliminaries}
\label{subsec-preliminary}

We first introduce some notations. For a subset $C\subset \R^n$, let $C^0$ and $\overline{C}$ denote the interior and closure of $C$, respectively.
For a vector $x=(x_1, \ldots, x_n)^\top\in R^n$, denote
$\|x\|_\infty = \max_{1\le i\le n} |x_i|$ and $\|x\|_1=\sum_i |x_i|$ by the $\ell_\infty$- and $\ell_1$-norm of $x$, respectively.
Let $B(x, \epsilon)=\{y: \| x-y\|_\infty \le \epsilon\}$ be an $\epsilon$-neighborhood of $x$.
For an $n\times n$ matrix $J=(J_{ij})$, let $\|J\|_\infty$ denote the matrix norm induced by the $\ell_\infty$-norm on vectors in $\R^n$, i.e.,
\[
\|J\|_\infty = \max_{x\neq 0} \frac{ \|Jx\|_\infty }{\|x\|_\infty}
=\max_{1\le i\le n}\sum_{j=1}^n |J_{ij}|,
\]
and $\|J\|$ be a general matrix norm.
Define the maximum absolute entry-wise norm: $\|J\|_{\max}=\max_{i,j}|J_{ij}|$.
We use the superscript ``*" to denote the true parameter under which the data are generated.
When there is no ambiguity, we omit the superscript ``*".

Recall the definitions of $\pi_{ij}$ and $\mu_{ij}(\theta, \gamma)$ in \eqref{eq-definition-nota}.
To explore asymptotic properties of the moment estimators,
we define a set of functions
\begin{equation}
 \label{definition-F-Q}
\begin{aligned}
F_i(\theta,\gamma)=\sum\limits_{k=1}^{n}\mu(\pi_{ik})-d_i,\quad i=1,\ldots,m,
\\
F_{m+j}(\theta,\gamma)=\sum\limits_{k=1}^{m}\mu(\pi_{kj})-b_j,\quad j=1,\ldots,n,
\\
Q(\theta, \gamma)= \sum_{i=1}^{m} \sum_{j=1}^{n}z_{ij} ( \mu_{ij}(\theta, \gamma) - x_{ij} ).
\end{aligned}
\end{equation}
Further, define
\begin{equation}
\label{eq-defini-F}
F(\theta, \gamma)=(F_1(\theta, \gamma), \ldots, F_n(\theta, \gamma),F_{n+1}(\theta, \gamma), \ldots,F_{m+n-1}(\theta, \gamma) )^\top.
\end{equation}
For any given $\gamma$, we use $F_{\gamma,i}(\theta)$ to denote $F_i(\theta, \gamma)$, and denote $F_\gamma(\theta)=(F_{\gamma,1}(\theta), \ldots, F_{\gamma,m+n-1}(\theta))^\top$.
Denote $\widehat\theta_\gamma$ by the solution to $F_\gamma(\theta)=0$.
We define a profiled function for $Q(\theta, \gamma)$:
\begin{equation}
\label{definition-Qc}
Q_c(\gamma)= \sum_{i=1}^{m} \sum_{j=1}^{n} z_{ij} ( \mu_{ij}(\widehat\theta_\gamma, \gamma) - x_{ij} ).
\end{equation}
According to the definitions, we have the following equations:
\begin{equation}\label{equation:FQ}
F(\widehat{\theta}, \widehat{\gamma})=0,~~F_\gamma(\widehat{\theta}_\gamma)=0,~~
Q(\widehat{\theta}, \widehat{\gamma})=0,~~Q_c(\widehat{\gamma})=0.
\end{equation}

We define a matrix class $ {\cal L}_{m,n}(M_0,M_1)$ with $ M_1\ge M_0 >0$.
For an $(m+n-1)\times (m+n-1)$ matrix $V=(v_{i,j})$,
we say it belongs to the matrix class $ {\cal L}_{m,n}(M_0,M_1)$ if the following conditions hold:
\begin{equation}\label{eq:Lmatrix}
\begin{split}
&M_0\leq v_{i,i}-\sum_{j=m+1}^{m+n-1}v_{i,j}\leq M_1,~~i=1,\dots,m, \\ 
&v_{i,j}=0,~~i,j=1,\dots,m,i\neq j,\\
&v_{i,j}=0,~~i,j=m+1,\dots,m+n-1,i\neq j,\\
&M_0\leq v_{i,j}=v_{j,i}\leq M_1,~~i=1,\dots,m,j=m+1,\dots,m+n-1,\\
&v_{i,i}=\sum_{k=1}^{m}v_{k,i}=\sum_{k=1}^{m}v_{i,k},~~i=m+1,\dots,m+n-1.
\end{split}
\end{equation}
It is easy to see that  $V$ is a diagonally dominant non-negative symmetric matrix.
Therefore, $V$ is positive definite.
When $\mu^\prime(x)>0$, the Jacobian matrix $F'_\gamma(\theta)$ belongs to this matrix class.
For convenience,  define
\begin{equation}
\begin{array}{rcl}
v_{m+n,i}=v_{i,m+n}&:=&v_{i,i}-\sum_{j=1;j\neq i}^{m+n-1}v_{i,j}, \quad i=1,\dots,m+n-1,
\\
v_{m+n,m+n}&:=&\sum_{i=1}^{m+n-1}v_{m+n,i}.
\end{array}
\end{equation}
In general, the inverse matrix $V^{-1}$ does not have a closed form.
\cite{fan2022asymptotic} proposed  a simple matrix  $S=(s_{i,j})$ to approximate $V^{-1}$, where
\begin{equation}\label{eq:Smatrix}
s_{i,j}=
\left\{
\begin{array}{ll}
\frac{\delta_{i,j}}{v_{i,i}}+\frac{1}{v_{m+n,m+n}}, & {i,j=1,\dots,m,} \\
-\frac{1}{v_{m+n,m+n}}, & {i=1,\dots,m,~j=m+1,\dots,m+n-1,}  \\
-\frac{1}{v_{m+n,m+n}}, & {i=m+1,\dots,m+n-1,j=1,\dots,m,}\\
\frac{\delta_{i,j}}{v_{i,i}}+\frac{1}{v_{m+n,m+n}}, & {i,j=m+1,\dots,m+n-1.}
  \end{array}
\right.
\end{equation}
In the above equation, $\delta_{ij}$ is a kronecker function, where it is equal to $1$ if $i=j$ and $0$ otherwise.

Let
\begin{equation}
\label{definition-H}
H(\theta, \gamma) := \frac{ \partial Q(\theta, \gamma) }{ \partial \gamma^T } - \frac{ \partial Q(\theta, \gamma) }{\partial \theta^T } \left[ \frac{\partial F(\theta, \gamma)}{\partial \theta^T } \right]^{-1}
\frac{\partial F(\theta, \gamma)}{\partial \gamma^T }.
\end{equation}
A direct calculation gives
\begin{align*}
    \frac{ \partial F_\gamma(\widehat{\theta}_\gamma) }{\partial \gamma^T}
    &=~
    \frac{ \partial F(\widehat{\theta}_\gamma, \gamma) }{\partial \theta^T}
    \frac{\partial \widehat{\theta}_\gamma }{\gamma^T} + \frac{\partial F(\widehat{\theta}_\gamma, \gamma)}{\partial \gamma^T} = 0,
    \\
    \frac{ \partial Q_c(\gamma)}{ \partial \gamma^T} &=~
    \frac{\partial Q(\widehat{\theta}_\gamma, \gamma)}{\partial \theta^T}
     \frac{\partial \widehat{\theta}_\gamma }{\gamma^T} + \frac{ \partial Q(\widehat{\theta}_\gamma, \gamma) }{ \partial \gamma^T}.
\end{align*}
As we can see, the Jacobian matrix $Q_c^\prime(\gamma)=\partial Q_c^\prime(\gamma)/\partial \gamma$ takes the following form:
\[
\frac{ \partial Q_c(\gamma) }{ \partial \gamma^T } = H( \widehat{\theta}_\gamma, \gamma).
\]
Asymptotic behaviours of $\widehat{\gamma}$ depends crucially on $Q_c^\prime(\gamma)$.
We assume that $H(\theta, \gamma)$ is positive definite. Otherwise, $\widehat{\gamma}$ will be ill-posed.
When $x_{ij}$ belongs to the exponential family distribution, $H(\theta, \gamma)$ is the Fisher information matrix of the concentrated likelihood function about $\gamma$ (e.g., page 126 of \citet{amemiya1985advanced}), where
the asymptotic variance of $\widehat{\gamma}$ is $H^{-1}(\theta, \gamma)$.
Thus, the positive definite assumption of $H(\theta, \gamma)$ is suitable.
 We define
\begin{equation}
    \kappa_{m,n}
    :=
    \sup_{\theta\in B(\theta^*, \epsilon_{m,n,1})} \| mn \cdot H^{-1}(\theta, \gamma^*)\|_\infty,
\end{equation}
where $B(\theta^*, \epsilon_{m,n,1}) = \{ \theta: \| \theta-\theta^*\|_\infty \le \epsilon_{m,n,1}\}$.

\subsection{Consistency}
\label{subsec-consis}

We present the conditions that guarantee consistency of the moment estimator.

\begin{condition}
\label{condition-0}
$m/n=O(1)$ and $m\ge n$.
\end{condition}

\begin{condition}
\label{assumption-1}
Suppose that there exist $b_{m,n,0}, b_{m,n,1}, b_{m,n,2}, b_{m,n,3}>0$ such that
\begin{subequations}
\begin{gather}
\min_{i,j} \mu^\prime(\pi_{ij}) \cdot \max_{i,j} \mu^\prime(\pi_{ij}) >0\label{ineq-mu-key0},\\
b_{m,n,0}\le \min_{i,j} |\mu^\prime(\pi_{ij})| \le \max_{i,j}|\mu^\prime(\pi_{ij})|\le b_{m,n,1}\label{ineq-mu-keya},  \\
\max_{i,j}|\mu^{\prime\prime}(\pi_{ij})| \le b_{m,n,2}\label{ineq-mu-keyb}, \\
\max_{i,j}|\mu^{\prime\prime\prime}(\pi_{ij})| \le b_{m,n,3}\label{ineq-mu-keyc}.
\end{gather}
\end{subequations}
hold for all $\theta \in B(\theta^*, \epsilon_{m,n,1}), \gamma\in B(\gamma^*, \epsilon_{m,n,2})$, where two positive numbers $\epsilon_{m,n,1}>0$ and $\epsilon_{m,n,2}>0$
satisfying $\epsilon_{m,n,1}=o(1)$ and $\epsilon_{m,n,2}=o(1)$.
\end{condition}

\begin{condition}\label{assumption-2}
Suppose $\max\limits_{i,j} \|z_{ij}\|_\infty\leq C_z$ for some universal constant $C_z$.
\end{condition}

\begin{condition}\label{assumption-3}
For each pair $(i,j)$ with $1\le i\le m, 1\le j\le n$, the distribution of $x_{ij}$ is sub-exponential with parameter $h_{ij}$, where we define $h_{\max}:=\max\limits_{i,j} h_{ij}$.
\end{condition}

Condition \ref{condition-0} requires that $m$ and $n$ in the same order.
It holds in real-world data sets. For instance,
there are $943$ users and $1682$ movies in the MovieLens 100K data set\footnote{Available at \url{https://grouplens.org/datasets/movielens/100k}}, where $m/n=0.56$.
As mentioned before, we set $m\ge n$ for convenience, which could be replaced with $n\le m$.
Condition \eqref{assumption-1} is mild, which holds in many commonly used probability distributions.
We illustrate it with an example. Consider the logistic distribution
for $f(\cdot)$, where $\mu(x)=e^x/(1+e^x)$.
 A direct calculation gives that
\[
\mu^\prime(x) = \frac{e^x}{ (1+e^x)^2 },~~  \mu^{\prime\prime}(x) = \frac{e^x(1-e^x)}{ (1+e^x)^3 },~~ \mu^{\prime\prime\prime}(x) = \frac{e^x(1-4e^x+e^{2x})}{ (1+e^x)^4 }.
\]
Since $y(1-y) \le 1/4$ when $y\in [0,1]$, and
\[
|\mu^{\prime\prime}(x)| \le \frac{e^x}{ (1+e^x)^2 } \times \left|\frac{(1-e^x)}{ (1+e^x) }\right|,~~~~
|\mu^{\prime\prime\prime}(x)| =
\frac{e^x}{ (1+e^x)^2 } \times \left| \left[ \frac{(1-e^x)^2}{ (1+e^x)^2 } - \frac{2e^x}{ (1+e^x)^2 }  \right]\right|
\]
we have
\begin{equation}\label{eq-mu-d-upper}
|\mu^\prime(x)| \le \frac{1}{4}, ~~ |\mu^{\prime\prime}(x)| \le \frac{1}{4},~~ |\mu^{\prime\prime\prime}(x)| \le \frac{1}{4}.
\end{equation}
Condition \ref{assumption-2} holds when covariates are bounded and is used in
\cite{graham2017econometric}, \cite{yan2019} and \cite{wang2023asymptotic}.
When the covariates are unbounded,
we can simply transform them to bounded variables by using
sigmoid or probit functions.
Condition \ref{assumption-3} restricts the range of probability density functions
$f(\cdot)$ to be sub-exponential, where is a widely applied distribution family,
covering many common distributions. It is also used in \cite{yan2016asymptoticsSIN}
and \cite{wang2023asymptotic}.
This condition can be replaced by any other conditions that guarantee
concentration inequalities for degrees.
Let $g=(d_{1},\cdots,d_{m}, b_{1},\cdots,b_{n-1})^{\top}$.
We now state the consistency result.

\begin{theorem}\label{theorem 2}
Let $V=\partial F(\theta^*, \gamma^*)/\partial \theta^\top$ and  $\sigma_{m,n}^2 = (m+n-1)^2 \| (V^{-1}-S) \mathrm{Cov}(g) (V^{-1}-S) \|_{\max}$, where $S$ is defined in \eqref{eq:Smatrix}.
Suppose Conditions \ref{condition-0}--\ref{assumption-3} hold, and
\begin{eqnarray}
\label{eq-theorema-ca}
 \frac{\kappa_{m,n}^2  b_{m,n,1}^{16} b_{m,n,2}}{b_{m,n,0}^{9}} (\frac{b_{m,n,2}h_{\max}^2}{b_{m,n,0}^{9}}+ \sigma_{m,n}) =o\left(\frac{m}{\log m}\right),
 \\
\label{eqt11*}
\frac{b_{m,n,1}^{4}b_{m,n,2}h_{\max}}{b_{m,n,0}^{6}}=o\left(\sqrt{\frac{m}{\log m}}\right),
\end{eqnarray}
Then the moment estimator $(\widehat{\beta},\widehat{\gamma})$ exists with high probability, and we further have
\begin{align*}
\| \widehat{\gamma} - \gamma^{*} \|_\infty &=
 O_p\left(\frac{\kappa_{m,n}  b_{m,n,1}^7 \log{m}}{m} (\frac{b_{m,n,2}h_{\max}^2}{b_{m,n,0}^{9}}+ \sigma_{m,n}) \right)=o_p(1),
\\
\| \widehat{\theta} - \theta^* \|_\infty &= O_p\left(\frac{b_{m,n,1}^{2}h_{\max}}{b_{m,n,0}^{3}}\cdot\sqrt{\frac{\log m}{m}}\right)=o_p(1).
\end{align*}
\end{theorem}

In exponential family distributions, $V=\mathrm{Cov}(g)$.
In this case, we have
\[
 \|(V^{-1}-S) V (V^{-1}-S)\|_{\max} \le \|V^{-1} - S\|_{\max} +
\|SVS-S\|_{\max}=O(\frac{b_{m,n,1}^2}{n^2b_{m,n,0}^3}),
\]
whose proof is in in the Supplementary material. It leads to
\begin{equation}
\label{sigma-simplify}
\sigma_{m,n}^2 = O(\frac{b_{m,n,1}^2}{b_{m,n,0}^3}).
\end{equation}
and condition \eqref{eq-theorema-ca} becomes
\begin{eqnarray}
\label{eq-theorema-ca-a}
 \frac{\kappa_{m,n}^2  b_{m,n,1}^{16} b_{m,n,2}}{b_{m,n,0}^{9}} (\frac{b_{m,n,2}h_{\max}^2}{b_{m,n,0}^{9}}+ \frac{b_{m,n,1}}{b_{m,n,0}^{3/2}}) =o\left(\frac{m}{\log m}\right).
\end{eqnarray}

\begin{remark}
If $b_{m,n,0}, b_{m,n,1}, b_{m,n,2}, \kappa_{m,n}, h_{\max}$ are constants, then the convergence rates of $\widehat{\beta}$ and $\widehat{\gamma}$ are
in orders of $O_p((\log n/n)^{1/2})$ and $O_p(\log n/n)$
The former is the same as in \cite{Chatterjee:Diaconis:Sly:2011}.
The latter has a faster convergence rate. This is due to that
we observe that $n(n-1)/2$ independent random variables
and the dimension of the parameter $\gamma$ is fixed.
\end{remark}

\subsection{Asymptotic normality}

We present the asymptotic normality of the moments estimators $(\widehat\theta,\widehat\gamma)$. Recall $g=(d_{1},\cdots,d_{m},\allowbreak b_{1},\cdots,b_{n-1})^{\top}$.
Let $g_{m+n}= \sum_{i=1}^m d_i -\sum_{j=1}^{n-1}b_j=b_n$, where $g$ denotes the observed degree sequence.

\begin{theorem}\label{theorem 3}
Assume the conditions in Theorem \ref{theorem 2} hold.
If
\begin{equation}\label{eqt3.1}
\frac{\kappa_{m,n}^2  b_{m,n,1}^{14} \log{m}}{m} (\frac{b_{m,n,2}h_{\max}^2}{b_{m,n,0}^{9}}+ \sigma_{m,n})^2
  =o\left( \frac{m}{\log m} \right),
\end{equation}
then for any fixed $i$,
\[
\widehat{\theta}_i - \theta^*_i = [S(g - \E g)]_i +  O_p\left(\frac{\kappa_{m,n}  b_{m,n,1}^{10} \log{m}}{b_{m,n,0}^3m} (\frac{b_{m,n,2}h_{\max}^2}{b_{m,n,0}^{9}}+ \sigma_{m,n}) \right),
\]
where
$$
[S(g - \E g)]_i=\begin{cases}
\frac{g_i- \E g_i}{v_{ii}}+\frac{g_{m+n}-\E g_{m+n}}{v_{m+n,m+n}}, &i=1,\ldots,m,\\
\frac{g_i- \E g_i}{v_{ii}}-\frac{g_{m+n}-\E g_{m+n}}{v_{m+n,m+n}}, &i=m+1,\ldots,m+n-1.
\end{cases}
$$
\end{theorem}

Define $U = (u_{i,j}) = \mathrm{Cov}(g)$, and
\begin{equation}\label{eq.condition1}
0<\eta_{m,n,4}:= \min_{i,j}\mathrm{Var}(x_{i,j})\leq \eta_{m,n,5}:=\max_{i,j} \mathrm{Var}(x_{i,j}).
\end{equation}
Because the elements of $X=(x_{ij})$ are independent, we have
$u_{i,i'} = \mathrm{Cov}(\sum_{j'}x_{i,j'}, \sum_{j'}x_{i',j'}) = 0$ for
$1\leq i \neq i' \leq m$.
Similarly, $u_{j,j'}=0$ for $m+1\leq j\neq j'\leq m+n-1$, and
$u_{i,j} = \text{Cov}(\sum_{j'}x_{i,j'},\sum_{i'}x_{i',j}) = \mathrm{Var}(x_{i,j})$
for $1\le i\le m$ and $1\le j \le n$.  It is easy to check that $U \in \mathcal{L}_{m,n}(\eta_{m,n,4}, \eta_{m,n,5})$. Define $u_{m+n,m+n}=\text{Var}(b_{n})=\text{Var}(g_{m+n})$.
If $x_{i,j}$ is a sub-exponential random variable with parameter $h_{ij}$,
then
\[
 \E |x_{ij}|^{3}\leq (3h_{ij})^{3}.
\]
If $\eta_{m,n,4}^{3/2}/m^{1/2}\to 0$, then
$u^{-3/2}_{i,i}\sum^{n}_{j=1}\E|x_{i,j}|^{3}\to 0$.
This verifies Lyapunov's condition.
By Lyapunov's Central Limit Theorem \citep[e.g.][ page 362]{vershynin2018},
if $\eta_{m,n,4}^{3/2}/m^{1/2}\to 0$ and $\min_{i=1,\ldots,m+n} u_{ii} \to \infty$, then for any fixed $i$ and $j$,
we have
\[
\frac{ d_i - \E d_i }{ u_{ii}^{1/2} } \rightsquigarrow N(0,1),
\qquad
\frac{ b_j - \E b_j }{ u_{j+m, j+m}^{1/2} } \rightsquigarrow N(0,1),
\]
as $n\to\infty$. Here, $\rightsquigarrow$ denotes ``convergence in distribution."
Observe that
\[
\frac{ u_{ii} }{ v_{ii}^2 } \ge O( \frac{ \eta_{m,n,1} }{ m b_{m,n,0}^2 } ).
\]
If
\begin{equation}
\label{eq-eta-mn-1}
\frac{\kappa_{m,n}  b_{m,n,1}^{10} \log{m}}{b_{m,n,0}^3m} (\frac{b_{m,n,2}h_{m,n}^2}{b_{m,n,0}^{9}}+ \sigma_{m,n})
\times \frac{ m^{1/2} b_{m,n,0} }{ \eta_{m,n,1}^{1/2} } = o(1),
\end{equation}
then for any fixed $i\in\{1,\ldots,m\}$,
$$
\hat{\xi}_i:=\frac{\widehat{\alpha}_i - \alpha^*_i }{\sqrt{\frac{u_{ii}}{v_{ii}^2}+\frac{u_{m+n,m+n}}{v_{m+n,m+n}^2}}}=
\frac{\frac{g_i- \E g_i}{v_{ii}}+\frac{g_{m+n}-\E g_{m+n}}{v_{m+n,m+n}}}{\sqrt{\frac{u_{ii}}{v_{ii}^2}+\frac{u_{m+n,m+n}}{v_{m+n,m+n}^2}}}+o_p(1),
$$
and for any fixed $j\in \{1,\ldots,n-1\}$,
$$
\hat{\zeta}_j:=\frac{\widehat{\beta}_j - \beta^*_j }{\sqrt{\frac{u_{m+j,m+j}}{v_{m+j,m+j}^2}+\frac{u_{m+n,m+n}}{v_{m+n,m+n}^2}}}\frac{\frac{b_j- \E b_j}{v_{jj}}+\frac{g_{m+n}-\E g_{m+n}}{v_{m+n,m+n}}}{\sqrt{\frac{u_{j+m,j+m}}{v_{j+m,j+m}^2}+\frac{u_{m+n,m+n}}{v_{m+n,m+n}^2}}}
+o_p(1).
$$
Then, we have the following corollary.

\begin{corollary}
\label{corollary1}
Suppose the conditions in Theorem \ref{theorem 3} holds.
If \eqref{eq-eta-mn-1} holds, then
for fixed $r$ and $s$, the vector
$(\hat{\xi}_1, \ldots, \hat{\xi}_r, \hat{\zeta}_1, \ldots, \hat{\zeta}_s)$
converges in distribution to the $r+s$-dimensional multivariate standard normal distribution.
\end{corollary}

The above corollary establishes the asymptotic normal distribution of $\widehat\theta$. We now state the asymptotic normality of $\widehat\gamma$.
For $i=1,\ldots,m$ and $j=1,\ldots,n$, define
$T_{ij} \in\mathbb R^{m+n-1}$ by a column vector with
 all entries zero except its $i$th and $(m+j)$th element being equal to $1$.
Let $e_k\in\mathbb R^{m+n-1}$, where the $k$th element is equal to 1 and all other elements are $0$.
 Define
\begin{equation}
\label{definition-sij}
\begin{array}{c}
V(\theta, \gamma)=\frac{ \partial F(\theta, \gamma) }{ \partial \theta^T }, ~~
V_{Q\theta}(\theta, \gamma) = \frac{ \partial Q(\theta, \gamma) }{ \partial \theta^T}, \\
s_{ij}(\theta, \gamma) = (x_{ij}-\E x_{ij}) ( z_{ij} + V_{Q\theta}(\theta, \gamma) [V(\theta,\gamma)]^{-1} T_{ij}).
\end{array}
\end{equation}
When evaluating $H(\theta,\gamma)$, $Q(\theta, \gamma)$, $V(\theta, \gamma)$ and $V_{Q\theta}(\theta, \gamma)$ at their true values
$(\theta^*, \gamma^*)$, we omit the arguments $\theta^*, \gamma^*$, i.e., $V=V(\theta^*, \gamma^*)$.
Let $N=mn$.  Also define
\[
\bar{H}=   \frac{1}{N} H( \theta^*, \gamma^*),
\]
where we recall the definition of $H(\theta, \gamma)$ from \eqref{definition-H}.
We now state the asymptotic representation of $\widehat{\gamma}$.

\begin{theorem}
\label{theorem 4}
Assume the conditions in Theorem \ref{theorem 2} hold.
If
\[
\frac{b_{m,n,3}b_{m,n,1}^6h_{\max}^3}{b_{m,n,0}^9} = o( \frac{m^{1/2}}{(\log m)^{3/2}}),
\]
then we have
\[
\sqrt{N}(\widehat{\gamma}- \gamma^*) = -\bar{H}^{-1} B_* + \bar{H}^{-1} \times \frac{1}{\sqrt{N}}   \sum_{i=1}^{m}\sum_{j=1}^{n}
s_{ij} (\theta^*, \gamma^*) + o_p(1).
\]
where
\begin{equation}\label{defintion-Bias}
B_*=\lim_{m\to\infty} \frac{1}{2\sqrt{N}} \sum_{k=1}^{m+n-1} \left[\frac{ \partial^2 Q(\theta^*, \gamma^*) }{ \partial \theta_k \partial \theta^\top}
V^{-1} U V^{-1} e_k \right].
\end{equation}
\end{theorem}

When $f(\cdot)$ belongs to the exponential family, it is easy to compute that $V=U$, where $\partial F(\theta^*, \gamma^*)/\partial \theta=\mathrm{Var}(g)$, In this case, $B_*$ simplifies to:
\begin{equation}\label{BB}
B_*=\lim_{m\to\infty} \frac{1}{2\sqrt{N}} \left[ \sum_{k=1}^{m}\left(  \frac{\sum_{j=1}^{n}z_{kj}\mu''(\pi_{kj})}{\sum_{j=1}^{n}\mu'(\pi_{kj})} \right)+
\sum_{k=m+1}^{m+n}\frac{\sum_{i=1}^{m}z_{i,k-m}\mu''(\pi_{i,k-m})}{\sum_{i}^{m}\mu'(\pi_{i,k-m})} \right].
\end{equation}

Let $\tilde{z}_{ij}=  z_{ij} + V_{Q\theta}V^{-1} T_{ij}$.
Recall that $\max\limits_{i,j} \|z_{ij}\|_\infty=O(1)$.
By calculations, we have $\|V_{Q\theta}\|_{\max}= O(mb_{m,n,1})$.
By setting $V^{-1}=S+W$, we have
\begin{eqnarray*}
\|V_{Q\theta}V^{-1}T_{ij}\|_{\max}
&\leq& 2\|V_{Q\theta}V^{-1}\|_{\max} =    2\|V_{Q\theta}(S+W)\|_{\max} \\
&\leq& 2(\|V_{Q\theta}S\|_{\max}+\|V_{Q\theta}W\|_{\max} )\\
&\leq& 2\left(\max_i\frac{\|V_{Q\theta}\|_{\max}}{v_{ii}}+\frac{1}{v_{m+n,m+n}}\max_{l=1,\ldots p}|\sum_{i=1}^m z_{inl}\mu^\prime(\pi_{in})|
\right.
\\
&&\left.              +(m+n-1)\|V_{Q\theta}\|_{\max}\|W\|_{\max}\right)\\
&=&O(\frac{b_{m,n,1}^3}{b_{m,n,0}^3}).
\end{eqnarray*}
So, we have $\|\tilde{z}_{ij}\|_{\max}=O(\frac{b_{m,n,1}^3}{b_{m,n,0}^3})$.
This shows that $\tilde{z}_{ij}$ is bounded.
For any nonzero vector $c=(c_1, \ldots c_p)^\top$, when $m,n \to \infty$, we have
\begin{equation}\label{eqt4-lyyaponu-condition}
\frac{ \sum\limits_{i=1}^m\sum\limits_{j=1}^n |c^\top \tilde{z}_{ij} |^3 \mathrm{E}|x_{ij}-\mathrm{E}x_{ij}|^3  }
{ [\sum\limits_{i=1}^m\sum\limits_{j=1}^n (c^\top \tilde{z}_{ij} )^2 \mathrm{Var}(x_{ij})]^{\frac{3}{2}}  } \to 0,
\end{equation}
then
$(\Sigma)^{-1/2}  \sum\limits_{i=1}^m\sum\limits_{j=1}^n s_{ij} (\theta^*, \gamma^*)$
converges in distribution to the multivariate standard normal distribution, where
$\Sigma= \sum\limits_{i=1}^m\sum\limits_{j=1}^n \mathrm{Cov}(s_{ij} (\theta^*, \gamma^*)) $.
We have the following corollary.

\begin{corollary}\label{corollary2}
Assume the conditions in Theorem \ref{theorem 4} and \eqref{eqt4-lyyaponu-condition} hold.
Then, we have
\begin{equation}
    \sqrt{N}(\widehat{\gamma} - \gamma)
\rightsquigarrow
    N\Big(
        -\bar{H}^{-1}B_*,
        \frac{1}{N}\bar{H}^{-1} \Sigma (\bar{H}^{-1})^\top
    \Big).
\end{equation}
\end{corollary}

\section{Applications}
\label{section:app}

We illustrate the theoretical result by using
two commonly used distributions for $f(\cdot)$ in model \eqref{model}.
the logistic distribution and Poisson distribution, where
the MLE is the same as the moment estimator.

\subsection{The logistic model}
\label{section-logistic}

When $f(\cdot)$ takes the logistic distribution, we have
\[
\P(a_{ij}=1) = \frac{  e^{\alpha_i + \beta_j + z_{ij}^\top \gamma  }}{ 1 + e^{ \alpha_i + \beta_j + z_{ij}^\top \gamma } }.
\]
We call it the covariate-adjusted bipartite $\beta$-model hereafter.

The numbers involved with the conditions are as follows.
Because $x_{ij}$'s are Bernoulli random variables, they are sub-exponential with $h_{ij}=1$.
The numbers $b_{m,n,0}, b_{m,n,1}, b_{m,n,2}$ and $b_{m,n,3}$ defined in Condition \eqref{assumption-2} are all $1/4$
as shown in \eqref{eq-mu-d-upper}.
The conditions \eqref{eq-theorema-ca} and \eqref{eqt11*} in Theorem \ref{theorem 2} becomes that
\begin{equation}\label{app-co-a}
\kappa_{m,n}/b_{m,n,0}^9  = o\left( \sqrt{\frac{n}{\log n}} \right),
\end{equation}
where $b_{m,n,0} = O(e^{ -2\|\theta^*\|_\infty - \|\gamma^*\|_\infty })$.
By Theorem \ref{theorem 2}, we have the following corollary.

\begin{corollary}
If \eqref{app-co-a} holds, then
\[
\|\widehat{\gamma}-\gamma^*\|_\infty = O_p\left( \frac{\kappa_{m,n}\log n}{ nb_{m,n,0}^9} \right),
\quad
\|\widehat{\beta} - \beta^*\|_\infty = O_p\left( \frac{1}{b_{m,n,0}^3}  \sqrt{\frac{\log n}{n}} \right).
\]
\end{corollary}

We discuss the condition and convergence rates related to the network density.
The expected network density is
\[
\rho_{m,n} :=\frac{1}{mn} \sum_{i,j} \E x_{ij} =
\frac{1}{N} \sum_{i,j}
\frac{ e^{\alpha_i + \beta_j + z_{ij}^\top \gamma } }
{ 1 + e^{\alpha_i + \beta_j + z_{ij}^\top \gamma } }.
\]
We consider one-dimensional covariate $z_{ij}$ for illustration.
In this case, by using $S$ in \eqref{eq:Smatrix} to approximate $V^{-1}$,
$H(\beta,\gamma)$ defined in \eqref{definition-H} can be written as
 \[
 H(\beta^*,\gamma^*)=\sum_{i,j} z_{ij}^2 \mu^\prime( \pi_{ij} ) -
 \sum_{i=1}^m \frac{1}{v_{ii}} \left(\sum_{j=1}^n z_{ij}\mu^\prime(\pi_{ij})\right)^2
 -\sum_{j=1}^n \frac{1}{v_{j+m,j+m}} \left(\sum_{i=1}^m z_{ij}\mu^\prime(\pi_{ij})\right)^2,
 \]
 such that $\kappa_{m,n}$ is approximately the inverse of $n^{-2}H(\beta^*, \gamma^*)$.
 It  depends on the covariates, the configuration of parameters, and the derivative of the mean function $\mu(\cdot)$.
 Therefore,
it is not possible to express $\kappa_{m,n}$ and $b_{m,n,0}$ as a function
with only  $\rho_n$ as its argument.
We consider one special case that $\theta_1 = \cdots = \theta_{m+n-1} \leq c $ with $c$ as a fixed constant, and assume that $z_{ij}$ is independently drawn from
the standard normality.
In this case, by large sample theory, we have
\[
\frac{1}{mn}\sum_{i,j} z_{ij}^2 \mu^\prime( \pi_{ij} ) \stackrel{p.}{\to}
\frac{e^{2\alpha_1}}{(1+e^{2\alpha_1})^2}, \qquad
\frac{1}{mn}\sum_{i,j} z_{ij}\mu^\prime(\pi_{ij}) \stackrel{p.}{\to} 0,
\]
such that $\kappa_n \asymp 1/\rho_n$, where $a_n \asymp b_n$ means
$c_1 a_n \le b_n \le c_2 a_n$ with two constants $c_1$ and $c_2$ for sufficiently large $n$.
Further, $b_{m,n,0} \asymp   \rho_n $.
Then the condition \eqref{eq-theorema-ca} becomes
\[
\frac{\rho_n}{ (\log n/n)^{1/20} }  \to \infty,
\]
and, the convergence rates are
\[
\| \widehat{\gamma} - \gamma^{*} \|_\infty =  O_p\left(
 \frac{  \log n }{ n \rho_{n}^{10}} \right  ),~~
\| \widehat{\beta} - \beta^* \|_\infty = O_p\left( \frac{ 1 }{\rho_{n}^3}\sqrt{\frac{\log n}{n}} \right).
\]
Here, estimation consistency requires a strong assumption $\rho_n \gg (n/\log n)^{1/8}$. It would be of interest to relax it.
The central limit theorems for $\widehat{\beta}$ and $\widehat{\gamma}$
directly follows from Corollaries \ref{corollary1} and \ref{corollary2}, respectively.

\subsection{The Poisson model}

We consider nonnegative integer weighth in Example \ref{example-b} i.e., $x_{ij}\in \{0, 1, \ldots\}$ and use the Poisson distribution for the weight $x_{ij}$, where
\[
\P(a_{ij}=k) = \frac{ \lambda_{ij}^k }{k!} e^{-\lambda_{ij}},
\]
and $\lambda_{ij}=  e^{z_{ij}^\top \gamma + \alpha_i + \beta_j }$.
Because it is an exponential family distribution, the MLE is identical to the moment estimator.
The expectation of $x_{ij}$ is $\lambda_{ij}=  e^{z_{ij}^\top \gamma + \alpha_i + \beta_j }$.
In this case, $\mu(x)=e^x$.
Define
\[
q_{m,n} := \sup_{\theta \in B_\infty(\theta^*, \epsilon_{m,n,1}), \gamma\in B_\infty(\gamma^*, \epsilon_{m,n,2}) }\max_{i,j} | \alpha_i + \beta_j + z_{ij}^\top \gamma |.
\]
So $b_{m,n,i}$'s ($i=0, \ldots, 3$) in inequalities \eqref{ineq-mu-keya}, \eqref{ineq-mu-keyb} and \eqref{ineq-mu-keyc} are
\[
b_{m,n,0} = e^{-q_{m,n}}, ~~ b_{m,n,1}= e^{q_{m,n}}, ~~ b_{m,n,2} = e^{q_{m,n}}, ~~
b_{m,n,3} = e^{q_{m,n}}.
\]
A Poisson distribution with parameter $\lambda$ is sub-exponential with parameter $c\lambda$
\citep[e.g.][example 4.6]{zhang2020concentration},
where $c$ is an absolute constant;
see Example 4.6 in \citet{zhang2020concentration}. Thus, $h_{ij}$ in Condition \ref{assumption-3} is $ce^{2q_{m,n}}$.
By Theorem \ref{theorem 2}, we have the following corollary.

\begin{corollary}
\label{coro-b}
Under Condition \ref{condition-0},
if $\kappa_{m,n} e^{37q_{m,n}} = o( (n/\log n)^{1/2})$, then
then
\[
\|\widehat{\gamma}-\gamma^*\|_\infty =  O_p(\frac{ \kappa_n e^{17 q_{m,n}} \log n}{n})=
o_p(1),~~~ \|\widehat{\beta} - \beta^*\|_\infty = O_p( \frac{ e^{5q_{m,n}}(\log n)^{1/2}}{n^{1/2}})=o_p(1).
\]
\end{corollary}

Note that $x_{ij}$'s $1\le i\le m, 1\le j \le n$ are
 independent Poisson random variables.
Since $v_{ij} = \E x_{ij} = \lambda_{ij}$, we have
\[
e^{-q_{m,n} } \le v_{ij}= e^{ \beta_i + \beta_j + z_{ij}^\top \gamma }
\le e^{q_{m,n}},~~1\le i< j\le n.
 \]
By using the Stein-Chen identity [\citet{Stein1972, chen1975}] for the Poisson distribution, it is easy to verify that
\begin{equation}\label{eqn:poisson:exp}
\E (a_{ij}^3) = \lambda_{ij}^3 + 3\lambda_{ij}^2 + \lambda_{ij}.
\end{equation}
It follows
\[
\frac{\sum_{j\neq i} \E (x_{ij}^3) }{ v_{ii}^{3/2} } \le \frac{ (n-1)e^{q_{m,n} } }{ (n-1)^{3/2} e^{-q_{m,n} } }
= O( \frac{ e^{4q_{m,n} } }{ n^{1/2} }).
\]
If $e^{4q_{m,n}}  = o( n^{1/2} )$, then the above expression goes to zero.
For any nonzero vector $c=(c_1, \ldots c_p)^T$, if
\begin{equation}\label{eqn:lemma5:a}
\frac{ \sum_{j<i } (c^\top \tilde{z}_{ij} )^3 \lambda_{ij}^3 }{ [\sum_{j<i } (c^\top \tilde{z}_{ij} )^2 \lambda_{ij} ]^{3/2} } = o(1),
\end{equation}
This verifies the condition \eqref{eqt4-lyyaponu-condition}.
Consequently, by Corollaries \ref{corollary1} and \ref{corollary2},
for fixed $r$ and $s$, the vector
$( \hat{\alpha}_1 - \alpha_1^*, \ldots,  \hat{\alpha}_r - \alpha_r^*,\hat{\beta}_1 - \beta_1^*, \ldots, \hat{\beta}_r - \beta_r^*)$ asymptotically follows a $(r+s)$-dimensional
normal distribution with covariance matrix $S_{(1,\ldots,r, m+1, \ldots, m+s),
(1,\ldots,r, m+1, \ldots, m+s)}$, where $S_{ \Omega_1, \Omega_2}$ denotes
the sub-matrix of $S$ on the row indexing $\Omega_1$ and column indexing $\Omega_2$.
In addition, $(mn)^{1/2} \overline{\Sigma}^{-1/2}(\hat{\gamma}-\gamma^*)$ converges in distribution to multivariate normal distribution with mean $\overline{\Sigma}^{-1/2}\bar{H}^{-1}B_*$  and covariance $I_p$,
where $I_p$ is the identity matrix, where $\bar{\Sigma}= N^{-1}\bar{H}^{-1} \tilde{\Sigma} \bar{H}^{-1}$.

\section{Numerical Analysis}
\label{sec-numerical}

\subsection{Simulation}
\label{subsec-simu}
To carry out simulations, we need to specify the distribution of $f(\cdot)$ in model \eqref{model}.
We consider the logistic distribution in section \ref{section-logistic}
and evaluate the asymptotic results of the estimator, where
the MLE is identical to the moment estimator.
We set the true parameter values in a linear form.
For $i=1, \ldots, m$, set $\alpha_{i}^* = (m-i)L/(m-1)$.
Similarly, for $j=1, \ldots, n$, set $\beta_{j}^* = (n-j)L/(n-1)$.
Here, $L$ is used to set different asymptotic regime, i.e.,
controlling the network density by changing $L$.
We consider four values of $L$: $L\in \{-0.2\log m, 0, 0.2\log m, 0.4\log m \}$.

The edge covariate between nodes $i$ and $j$ is generated according to the formula  $z_{ij}=(x_{i1} \times x_{j1}, x_{i2} \times x_{j2})^\top$, where $x_{i1}$ takes values $1$ or $-1$ with probabilities $0.3$ and $0.7$, $x_{j1}$ takes values $1$ or $-1$ with probabilities $0.6$ and $0.4$, $x_{i2}$ and $x_{j2}$ take values $1$ or $-1$ with equal probability.
All covariates are generated independently.
The regression coefficient of covariates is set to $\gamma^*=(0.5, 1)^\top$.
This setting posits positive effects of covariates, where
the edges are easier to occur between actors and events with matching attributes.

We first evaluate the average absolute error of the estimator.
We conducted simulations under two settings: $m=100,n=100$ and $m=300,n=100$.
Each simulation was repeated $5000$ times.
The simulation results are reported in Table \ref{tab1}, where
we only select several values for $\alpha$ and $\beta$.
From this table, we can see that
the estimation error becomes smaller as $m$ increases when $n$ is fixed.
The error of $\hat{\gamma}$ is significantly smaller than those of $\hat{\alpha}$ and $\hat{\beta}$.
These observations confirm the consistency result in Theorem \ref{theorem 2}.

\begin{table}[!htpb]
\centering
\caption{Mean $|\hat{\alpha}-\alpha|$ / Mean $|\hat{\beta}-\beta|$ / Mean $|\hat{\gamma}-\gamma|$.}
\label{tab1}
\small 
\renewcommand{\arraystretch}{0.6} 
\begin{tabular}{l c c c c c}
\toprule
$value$      & $L =-0.2\log m$ & $L = 0$  & $L=0.2\log m$ &  $L =0.4\log m$\\
\midrule
\multicolumn{5}{c}{$m = 100$, $n=100$} \\
\midrule
$\vert\hat{\alpha}_{1}-\alpha_{1}\vert$         &$ 0.281 $&$ 0.261 $&$ 0.283 $&$ 0.283 $ \\
$\vert\hat{\alpha}_{m/2}-\alpha_{m/2}\vert$     &$ 0.267 $&$ 0.260 $&$ 0.270 $&$ 0.270 $ \\
$\vert\hat{\alpha}_{m}-\alpha_{m}\vert$         &$ 0.269 $&$ 0.264 $&$ 0.264 $&$ 0.264 $ \\
$\vert\hat{\beta}_{1}-\beta_{1}\vert$           &$ 0.282 $&$ 0.265 $&$ 0.292 $&$ 0.292 $ \\
$\vert\hat{\beta}_{n/2}-\beta_{n/2}\vert$       &$ 0.282 $&$ 0.265 $&$ 0.277 $&$ 0.277 $\\
$\vert\hat{\beta}_{n-1}-\beta_{n-1}\vert$       &$ 0.273 $&$ 0.260 $&$ 0.262 $&$ 0.262 $ \\
$\vert\hat{\gamma}_{1}-\gamma_{1}\vert$         &$ 0.024 $&$ 0.023 $&$ 0.025 $&$ 0.025 $  \\
$\vert\hat{\gamma}_{2}-\gamma_{2}\vert$         &$ 0.029 $&$ 0.028 $&$ 0.029 $&$ 0.029 $ \\
\midrule
\multicolumn{5}{c}{$m = 300$, $n=100$} \\
\midrule
$\vert\hat{\alpha}_{1}-\alpha_{1}\vert$         &$ 0.265 $&$ 0.214 $&$ 0.256 $&$ 0.256 $ \\
$\vert\hat{\alpha}_{m/2}-\alpha_{m/2}\vert$     &$ 0.225 $&$ 0.218 $&$ 0.225 $&$ 0.225 $ \\
$\vert\hat{\alpha}_{m}-\alpha_{m}\vert$         &$ 0.218 $&$ 0.215 $&$ 0.220 $&$ 0.220 $ \\
$\vert\hat{\beta}_{1}-\beta_{1}\vert$           &$ 0.160 $&$ 0.150 $&$ 0.162 $&$ 0.162 $ \\
$\vert\hat{\beta}_{n/2}-\beta_{n/2}\vert$       &$ 0.155 $&$ 0.151 $&$ 0.154 $&$ 0.154 $\\
$\vert\hat{\beta}_{n-1}-\beta_{n-1}\vert$       &$ 0.156 $&$ 0.153 $&$ 0.157 $&$ 0.157 $ \\
$\vert\hat{\gamma}_{1}-\gamma_{1}\vert$         &$ 0.015 $&$ 0.014 $&$ 0.015 $&$ 0.015 $  \\
$\vert\hat{\gamma}_{2}-\gamma_{2}\vert$         &$ 0.019 $&$ 0.017 $&$ 0.018 $&$ 0.018 $ \\
\bottomrule
\end{tabular}
\end{table}

We now evaluate asymptotic distribution of the moment estimator.
According to Theorem~\ref{theorem 3}, for any $i,j=1, \ldots, m$, the normalized estimators
\begin{equation}
\label{test-st}
\widehat{\zeta}_i = \frac{\widehat{\alpha}_i - \alpha_i^*}{(1/\widehat{v}_{i,i} + 1/\widehat{v}_{m+n,m+n})^{1/2}}
\quad \text{and} \quad
\widehat{\xi}_{i,j} = \frac{\widehat{\alpha}_i - \widehat{\alpha}_{j} - (\alpha_i^* - \alpha_{j}^*)}{(1/\hat{v}_{i,i} + 1/\widehat{v}_{j,j})^{1/2}}
\end{equation}
converge in distribution to the standard normal distribution, where $\hat{v}_{i,i}$ is the estimator of $v_{i,i}$ obtained by replacing $(\theta^*, \gamma^*)$ with $(\widehat{\theta}, \widehat{\gamma})$. For the estimator $\widehat{\beta}$, there are parallel results.
However, we only show the asymptotic distribution of $\widehat{\alpha}$ for brevity.
Each simulation was repeated $5000$ times.

\begin{figure}
   \centering
    \includegraphics{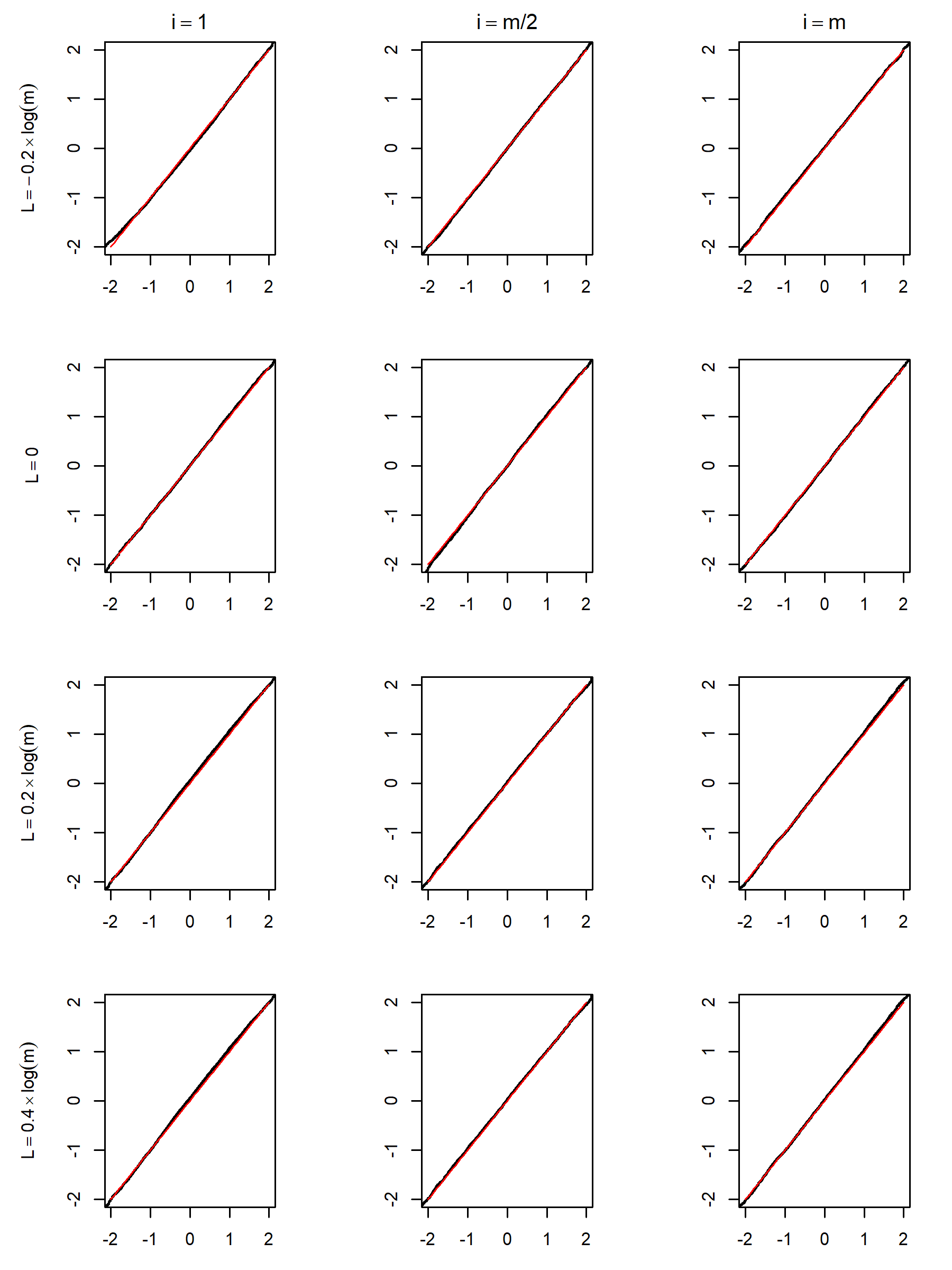}
    \caption{QQ plot of $\hat{\zeta}_i$ for $m=300$, $n=100$}
    \label{fig1}
\end{figure}

We employ QQ plots to assess the asymptotic normalities of $\widehat{\zeta}_i$
and $\widehat{\xi}_{i,j}$.
The plots are similar for $\widehat{\zeta}_i$
and $\widehat{\xi}_{i,j}$. We only show the figures for $\widehat{\zeta}_i$ in Figure \ref{fig1} to save space.
To compensate it, we report
the $95\%$ coverage frequency for $\alpha_i - \alpha_j$ in Table \ref{tab2}.
From Figure \ref{fig1}, we can see that
the sample quantiles are very close to the theoretical quantiles, where the QQ plots align with
the reference line $y=x$.
Table \ref{tab2} shows that
the simulated coverage frequencies are close to the target $95\%$ level.

\begin{table}[h]
\centering
\caption{Coverage probability ($\times 100\%$) / Confidence interval length for $\alpha_i - \alpha_j$}
\label{tab2}
\small 
\renewcommand{\arraystretch}{0.6} 
\begin{tabular}{ccccccc}
\hline
$(m,n)$      &  $(i,j)$ & $L=-0.2\log m$ & $L=0$ & $L=0.2\log m$ & $L=0.4\log m$ \\
\hline
(100,100)         &$(1,2)$      &$ 95.16 / 1.42 $&$ 95.2 / 1.28 $&$ 94.74 / 1.42 $&$ 94.74 / 1.42 $ \\
            &$(50,51) $         &$ 95.24 / 1.33 $&$ 94.78 / 1.28 $&$ 95.02 / 1.36 $&$ 95.02 / 1.36 $ \\
            &$(99,100)$         &$ 94.38 / 1.31 $&$ 94.60 / 1.28 $&$ 95.60 / 1.28 $&$ 95.60 / 1.28 $ \\
&&&&&&\\
(300,100)         &$(1,2)$      &$ 95.50 / 1.57 $&$ 94.48 / 1.27 $&$ 94.86 / 1.53 $&$ 94.86 / 1.53 $ \\
            &$(150,151) $       &$ 94.88 / 1.34 $&$ 94.72 / 1.27 $&$ 95.30 / 1.38 $&$ 95.30 / 1.38 $ \\
            &$(299,300) $       &$ 94.84 / 1.34 $&$ 94.32 / 1.27 $&$ 94.62 / 1.30 $&$ 94.62 / 1.30 $ \\
\hline
\end{tabular}
\end{table}

\begin{table}[h]
\centering
\caption{Coverage probability ($\times 100\%$) / Confidence interval length for $\gamma$}
\label{tab3}
\small 
\renewcommand{\arraystretch}{0.6} 
\begin{tabular}{ccccccc}
\hline
$(m,n)$     &   $\widehat{\gamma}$  & $L=-0.2\log m$ & $L=0$ & $L=0.2\log m$ & $L=0.4\log m$ \\
\hline
$(100,100)$   & $\hat{\gamma}_1$      &$ 94.82 / 0.11 $&$ 94.78 / 0.10 $&$ 94.20 / 1.36 $&$ 94.20 / 1.36 $ \\

        & $\hat{\gamma}_{bc, 1}$      &$ 94.50 / 0.11 $&$ 94.60 / 0.10 $&$ 94.50 / 1.36 $&$ 94.50 / 1.36 $ \\

        & $\hat{\gamma}_2$            &$ 94.78 / 0.11 $&$ 95.04 / 0.09 $&$ 94.44 / 1.32 $&$ 94.44 / 1.32 $ \\

        & $\hat{\gamma}_{bc, 2}$      &$ 94.74 / 0.11 $&$ 94.4 / 0.09 $&$ 95.42 / 1.32 $&$ 95.42 / 1.32 $\\

$(300,100)$   & $\hat{\gamma}_1$      &$ 92.64 / 0.07 $&$ 91.94 / 0.06 $&$ 92.22 / 0.07 $&$ 92.22 / 0.07 $ \\
        & $\hat{\gamma}_{bc, 1}$      &$ 95.08 / 0.07 $&$ 94.26 / 0.06 $&$ 94.92 / 0.07 $&$ 94.92 / 0.07 $ \\

        & $\hat{\gamma}_2$            &$ 84.82 / 0.07 $&$ 78.56 / 0.05 $&$ 82.94 / 0.06 $&$ 82.94 / 0.06 $ \\
        & $\hat{\gamma}_{bc, 2}$      &$ 93.82 / 0.07 $&$ 95.54 / 0.05 $&$ 93.22 / 0.06 $&$ 93.22 / 0.06 $\\

\hline
\end{tabular}
\end{table}

We now compare the uncorrected estimator $\hat{\gamma}$ with the bias-corrected estimator
$\hat{\gamma}_{\text{bc}} = \hat{\gamma} + \bar{H}^{-1}B_*$.
According to Corollary~\ref{corollary2}, the estimator $\hat{\gamma}$ exhibits a bias.
The results are reported in Table \ref{tab3}.
From this  table, we can see that (1) when $m=100$ and $n=100$,
the coverage frequencies are close to the target level $95\%$ for both $\hat{\gamma}$
and $\hat{\gamma}_{\text{bc}}$; (2) when $m=300$ and $n=100$,
the coverage frequencies are far lower than the target level $95\%$
while $\hat{\gamma}_{\text{bc}}$ achieves coverage frequencies close to $95\%$.

\subsection{A real Data example}
\label{subsec-real}

We use the covariate-adjusted bipartite $\beta$-model to analyze the
the MovieLens data \citep{harper2015movielens},
available from \url{https://grouplens.org/datasets/movielens/100k/}.
This data set contains $100,000$ ratings from $943$ users on $1682$ movies, along with user's demographics (e.g., age, gender) and movie's genres.
We construct a binary bipartite network, where an edge exists
between an user and a movie if and only if the user had ranked the movie.
We chose those users and movies having degrees over $40$ and constructed the
reduced subgraph for subsequent analysis.
The resulting subgraph has $645$ users and $708$ movies.
We generates two covariate variable for each pair of users and movies.
For the first covariate variable $z_{ij,1}$ between user $i$ and movie $j$,
we formulate it by matching user's sex and types of movies,
where we classify all the $18$ different types of movies
(excluding the unknown movie genre) into two groups ($G_{\mathrm{male}}$ and $G_{\mathrm{female}}$) according
the sex of users (male or female).
Specifically, $z_{ij,1}=1$ if user $i$ is male and movie $j$ genre belongs to $G_{\mathrm{male}}$ (or user $i$ is female and movie $j$ genre belongs to $G_{\mathrm{female}}$;
$z_{ij,1}=0$ otherwise.
For the second covariate variable $z_{ij,2}$, the definition is similar,
where all ages are categorized into three classes ($\le 18$, $(18, 55]$, $>55$)
and all types of movies are correspondingly classified into three groups.

\begin{figure}
\centering
\includegraphics[scale=0.7]{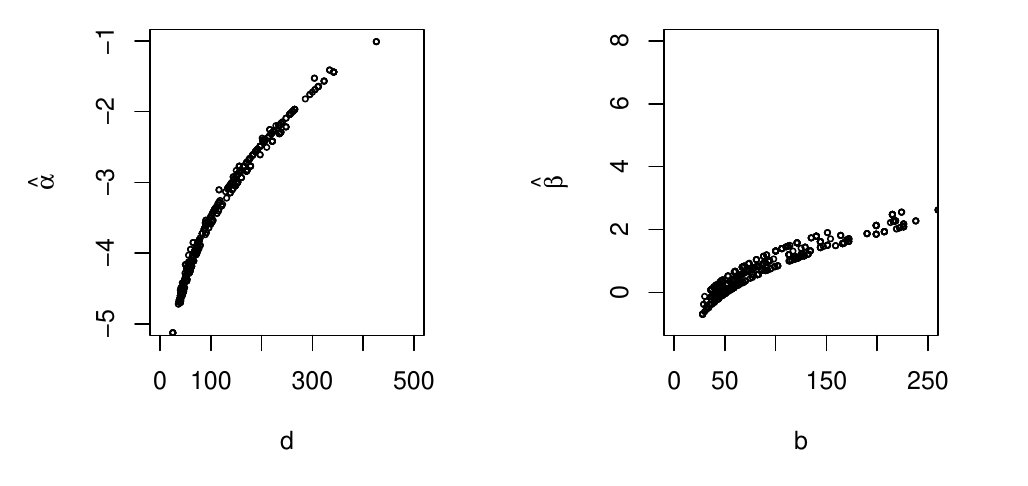}
\caption{Plots of the MLEs for the degree parameters.}
\label{fig-real-a}
\end{figure}

\begin{figure}
\centering
\includegraphics[scale=0.7]{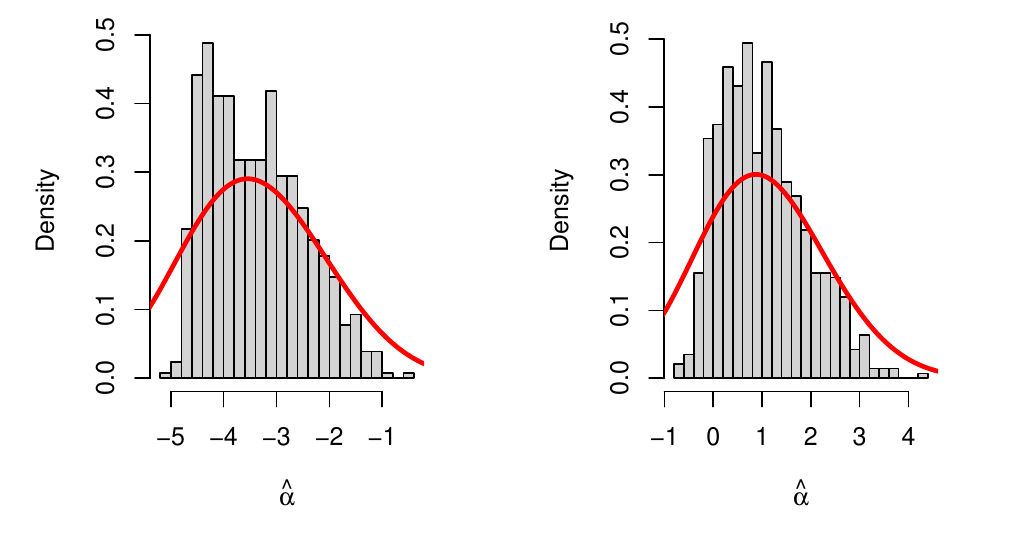}
\caption{Histogram density estimation of the fitted parameter estimates. }
\label{fig-hist}
\end{figure}

The estimated values of $\hat{\alpha}_i$ and $\hat{\beta}_j$
are shown in Figure \ref{fig-real-a}.
The left subfigure is about the degree $d_i$ of a user against the MLE
$\hat{\alpha}_i$,
while the right subfigure is about the degree $b_i$ of a movie  against the MLE
$\hat{\beta}_j$, where we arranged the degrees in a increasing manner.
This figure exhibits a large node heterogeneity phenomenon.
In addition, those nodes with high degrees have relatively large parameter estimates.
Figure \ref{fig-hist} shows the histogram density estimation of the fitted parameter estimates,
where the density curves for $\hat{\alpha}$ and $\hat{\beta}$ looks different.

We use some specific hypotheses as illustrating examples.
Consider the null hypothesis $H_{01}:
\alpha_1=\alpha_2$ and $H_{02}: \alpha_2=\alpha_3$.
We use the test statistic in \eqref{test-st} to calculate the $p$-values.
For $H_{01}$, the test statistic has a value $12.7$ with a p-value less than $10^{-3}$
and the standard error $0.17$.
For $H_{01}$, the test statistic has a value $0.82$ with a p-value $0.41$
and the standard error $0.22$.
This shows that it is significant for testing $\alpha_1=\alpha_2$
while it is not significant for testing $\alpha_2=\alpha_3$.
The fitted regression coefficients of covariates are
 $\hat{\gamma}_{1} = 0.359$ and $\hat{\gamma}_{bc,2} = 0.252$
with standard errors $0.013$ and $0.022$.
For testing $\gamma_1=0$ and $\gamma_2=0$, we obtain
both $p$-values less than $10^{-3}$, indicating
a significant covariate effect.

\section{Summary and discussion}
\label{section:sd}

We have developed the moment estimation in a general covariate-adjusted bipartite network model, which has a node-specific parameter for each node and
a fixed dimensional regression coefficient for covariates.
When both numbers of actors and events go to infinity,
we have established the consistency and asymptotic representation of the moment estimator under some conditions.
The conditions in \eqref{eq-theorema-ca} and \eqref{eqt11*}
imply that the network density can not be very small.
However, our simulations indicate that these conditions could be relaxed.
The asymptotic behavior of the moment estimator depends not only on the ranges of parameters,
but also on the configuration of all the parameters.
It would be of interest to see whether these conditions could be relaxed.

We do not consider the connections between actors or between events in bipartite networks.
If two actors are friends, then they tend to join in the same events likely.
Similarly, there may also be some connections between events with the same attributes,
which have influence on edge formations between actors and events.
This leads to that some events may have many common actors.
We do not consider these factors in this paper.
To address this issue, new models need to be developed, by taking into account
both one-mode and two-mode network features.
We would like to investigate this problem in the future.

\end{document}